\documentclass{amsart}
\usepackage{curvesls}
\usepackage{epic}
\usepackage{amsfonts}

\theoremstyle{plain}
\newtheorem{thm}{Theorem}[section]
\newtheorem{lem}[thm]{Lemma}

\newtheorem*{clm}{Claim}

\theoremstyle{definition}
\newtheorem*{defn}{Definition}

\newcommand{\Cplx}{{\bf{C}}}
\newcommand{\Real}{{\bf{R}}}
\newcommand{\Zed}{{\bf{Z}}}

\newcommand{\su}{{\mathfrak{su}}}
\newcommand{\lf}{{\ell f}}
\newcommand{\Hom}{{\mathrm{Hom}}}
\newcommand{\Local}{{\mathcal{L}}}
\newcommand{\pcolor}{{\mathbf{p}}}
\newcommand{\mcolor}{{\mathbf{m}}}

\newcommand{\suchthat}{\thinspace | \thinspace}

\newcommand{\stdspace}{\hskip 0.75em plus 0.15em \ignorespaces}
\newcommand{\ppar}{\par\goodbreak\vskip 8pt plus 3pt minus 3pt}
\newcommand{\co}{\colon\thinspace}
\newcommand{\rk}[1]{\ppar{\bf #1}\stdspace}

\newcommand{\lplus}
{\raisebox{-8pt}
  {
  \begin{picture}(20,20)
  \put(0,0){\vector(1,1){20}}
  \drawline(20,0)(12,8)
  \drawline(8,12)(0,20)
  \put(0,20){\vector(-1,1){0}}
  \end{picture}
  }
}
\newcommand{\lminus}
{\raisebox{-8pt}
  {
  \begin{picture}(20,20)
  \put(20,0){\vector(-1,1){20}}
  \drawline(0,0)(8,8)
  \drawline(12,12)(20,20)
  \put(20,20){\vector(1,1){0}}
  \end{picture}
  }
}
\newcommand{\lzero}
{\raisebox{-8pt}
  {
  \begin{picture}(20,20)
  \qbezier(0,0)(8,8)(8,10)
  \qbezier(8,10)(8,12)(0,20)
  \qbezier(20,0)(12,8)(12,10)
  \qbezier(12,10)(12,12)(20,20)
  \put(0,20){\vector(-1,1){0}}
  \put(20,20){\vector(1,1){0}}
  \end{picture}
  }
}
\newcommand{\eights}
{
  {
  \begin{picture}(120,60)
  \put(40,30){\circle*{4}}
  \qbezier(50,30)(30,10)(30,30)  
  \qbezier(50,30)(30,50)(30,30)  
  \qbezier(60,30)(20,-10)(20,30) 
  \qbezier(60,30)(20,70)(20,30)
  \qbezier(70,30)(10,-30)(10,30) 
  \qbezier(70,30)(10,90)(10,30)

  \put(80,30){\circle*{4}}
  \qbezier(70,30)(90,10)(90,30)  
  \qbezier(70,30)(90,50)(90,30)  
  \qbezier(60,30)(100,-10)(100,30) 
  \qbezier(60,30)(100,70)(100,30)
  \qbezier(50,30)(110,-30)(110,30) 
  \qbezier(50,30)(110,90)(110,30)

  \put(10,30){\vector(0,-1){0}}
  \put(20,30){\vector(0,-1){0}}
  \put(30,30){\vector(0,-1){0}}
  \put(0,25){$\gamma_3$}

  \put(90,30){\vector(0,-1){0}}
  \put(100,30){\vector(0,-1){0}}
  \put(110,30){\vector(0,-1){0}}
  \put(111,25){$\gamma_1$}

  \end{picture}
  }
}

\newcommand{\baab}
{
  {
  \begin{picture}(40,40)
  \put(10,0){\line(0,1){23}}
  \put(10,27){\line(0,1){13}}

  \put(20,0){\line(0,1){28}}
  \put(20,32){\line(0,1){8}}

  \qbezier(30,0)(30,5)(22,9)
  \drawline(18,11)(12,14)
  \qbezier(8,16)(0,20)(8,24)
  \put(8,24){\line(2,1){14}}
  \qbezier(22,31)(30,35)(30,40)

  \put(40,0){\line(0,1){40}}
  \end{picture}
  }
}

\newcommand{\bb}
{
  {
  \begin{picture}(30,40)
  \put(0,0){\line(0,1){40}}

  \put(10,0){\line(0,1){5}}
  \qbezier(10,5)(10,10)(13,13)
  \drawline(13,13)(17,17)
  \qbezier(17,17)(20,20)(17,23)
  \qbezier(13,27)(10,30)(10,35)
  \put(10,35){\line(0,1){5}}

  \put(20,0){\line(0,1){5}}
  \qbezier(20,5)(20,10)(17,13)
  \qbezier(13,17)(10,20)(13,23)
  \drawline(13,23)(17,27)
  \qbezier(17,27)(20,30)(20,35)
  \put(20,35){\line(0,1){5}}

  \put(30,0){\line(0,1){40}}
  \end{picture}
  }
}

\newcommand{\babCC}
{
  {
  \begin{picture}(50,60)

  \put(30,0){\line(0,1){5}}
  \qbezier(30,5)(30,10)(22,14)
  \drawline(18,16)(12,18)
  \qbezier(8,20)(0,24)(8,28)
  \put(8,28){\line(2,1){30}}
  \qbezier(42,45)(50,49)(42,53)
  \drawline(42,53)(38,55)
  \qbezier(38,55)(30,59)(30,65)
  \put(30,65){\line(0,1){5}}

  \put(10,0){\line(0,1){27}}
  \put(10,31){\line(0,1){39}}

  \put(20,0){\line(0,1){32}}
  \put(20,36){\line(0,1){34}}

  \put(40,0){\line(0,1){52}}
  \put(40,56){\line(0,1){14}}
  \end{picture}
  }
}

\newcommand{\fulltwist}
{
  \setlength{\unitlength}{0.7pt}
  {
  \begin{picture}(160,190)  

  \put(40,0){\line(0,1){15}}
  \qbezier(40,15)(40,20)(48,24)
  \put(48,24){\line(2,1){24}}
  \qbezier(72,36)(80,40)(80,45)
  \put(80,45){\line(0,1){25}}
  \qbezier(80,70)(80,75)(72,79)
  \put(68,81){\line(-2,1){36}}
  \qbezier(28,101)(20,105)(20,110)
  \put(20,110){\line(0,1){8}}
  \put(20,122){\line(0,1){6}}
  \put(20,132){\line(0,1){3}}
  \qbezier(20,135)(20,140)(28,144)
  \drawline(28,144)(32,146)
  \qbezier(32,146)(40,150)(40,155)
  \put(40,155){\line(0,1){5}}

  \put(50,0){\line(0,1){10}}
  \qbezier(50,10)(50,15)(58,19)
  \put(58,19){\line(2,1){24}}
  \qbezier(82,31)(90,35)(90,40)
  \put(90,40){\line(0,1){15}}
  \qbezier(90,55)(90,60)(82,64)
  \drawline(78,66)(72,69)
  \put(68,71){\line(-2,1){36}}
  \drawline(28,91)(18,96)
  \qbezier(18,96)(10,100)(10,105)
  \put(10,105){\line(0,1){8}}
  \put(10,117){\line(0,1){3}}
  \qbezier(10,120)(10,125)(18,129)
  \put(18,129){\line(2,1){24}}
  \qbezier(42,141)(50,145)(50,150)
  \put(50,150){\line(0,1){10}}

  \put(60,0){\line(0,1){5}}
  \qbezier(60,5)(60,10)(68,14)
  \put(68,14){\line(2,1){24}}
  \qbezier(92,26)(100,30)(100,35)
  \put(100,35){\line(0,1){5}}
  \qbezier(100,40)(100,45)(92,49)
  \drawline(88,51)(82,54)
  \drawline(78,56)(72,59)
  \put(68,61){\line(-2,1){36}}
  \put(28,81){\line(-2,1){20}}
  \qbezier(8,91)(0,95)(0,100)
  \put(0,100){\line(0,1){5}}
  \qbezier(0,105)(0,110)(8,114)
  \put(8,114){\line(2,1){44}}
  \qbezier(52,136)(60,140)(60,145)
  \put(60,145){\line(0,1){15}}

  \put(30,0){\line(0,1){123}}
  \put(30,127){\line(0,1){6}}
  \put(30,137){\line(0,1){6}}
  \put(30,147){\line(0,1){13}}
  
  \put(70,0){\line(0,1){13}}
  \put(70,17){\line(0,1){6}}
  \put(70,27){\line(0,1){6}}
  \put(70,37){\line(0,1){123}}

  \multiput(120,0)(10,0){5}{\line(0,1){160}}
  \end{picture}
  }
}
\newcommand{\bigx}
{
  {
  \begin{picture}(100,60)

  \put(0,0){\line(0,1){60}}

  \put(10,0){\line(0,1){10}}
  \qbezier(10,10)(10,15)(18,19)
  \put(18,19){\line(2,1){54}}
  \qbezier(72,46)(80,50)(80,55)
  \put(80,55){\line(0,1){5}}

  \put(20,0){\line(0,1){5}}
  \qbezier(20,5)(20,10)(28,14)
  \put(28,14){\line(2,1){54}}
  \qbezier(82,41)(90,45)(90,50)
  \put(90,50){\line(0,1){10}}

  \put(30,0){\line(0,1){13}}
  \put(30,17){\line(0,1){6}}
  \put(30,27){\line(0,1){33}}

  \put(70,0){\line(0,1){33}}
  \put(70,37){\line(0,1){6}}
  \put(70,47){\line(0,1){13}}

  \put(80,0){\line(0,1){5}}
  \qbezier(80,5)(80,10)(72,14)
  \drawline(68,16)(52,24)
  \drawline(48,26)(42,29)
  \drawline(38,31)(32,34)
  \drawline(28,36)(18,41)
  \qbezier(18,41)(10,45)(10,50)
  \put(10,50){\line(0,1){10}}

  \put(90,0){\line(0,1){10}}
  \qbezier(90,10)(90,15)(82,19)
  \drawline(82,19)(72,24)
  \drawline(68,26)(62,29)
  \drawline(58,31)(52,34)
  \drawline(48,36)(32,44)
  \qbezier(28,46)(20,50)(20,55)
  \put(20,55){\line(0,1){5}}
  
  \put(100,0){\line(0,1){60}}
  \end{picture}
  }
}
\newcommand{\transparent}
{
  {
  \begin{picture}(60,48)
  \put(12,24){\circle*{4}}
  \put(48,24){\circle*{4}}
  \put(12,24){\line(1,0){36}}

  \put(12,12){\line(1,0){36}}
  \put(48,24){\oval(24,24)[r]}
  \put(48,36){\line(-1,0){36}}
  \put(12,24){\oval(24,24)[l]}

  \put(24,0){\line(0,1){48}}
  \put(36,0){\line(0,1){48}}

  \put(13,39){$a_1$}
  \put(13,27){$a_2$}
  \put(13,15){$a_3$}
  \put(37,39){$b_1$}
  \put(37,27){$b_2$}
  \put(37,15){$b_3$}
  \end{picture}
  }
}
\newcommand{\barify}
{
  {
  \begin{picture}(50,100)

  \drawline(0,0)(0,100)(20,100)(20,20)(35,5)(45,5)
           (45,95)(35,95)(35,20)(15,0)(0,0)
  \drawline(5,5)(5,95)(15,95)(15,20)(35,0)(50,0)
           (50,100)(30,100)(30,20)(15,5)(5,5)
  \put(10,50){\circle*{4}}
  \put(40,50){\circle*{4}}
  \end{picture}
  }
}
\newcommand{\markov}
{
  {
  \begin{picture}(100,100)
  \put(0,0){\line(0,1){100}}
  \put(0,0){\line(1,0){50}}
  \put(0,100){\line(1,0){50}}
  \put(50,50){\oval(100,100)[r]}
  \put(20,50){\circle*{4}}
  \put(50,50){\circle*{4}}
  \put(80,50){\circle*{4}}

  \qbezier(50,0)(60,30)(60,50)
  \put(45,50){\oval(30,60)[t]}
  \put(50,50){\line(-1,0){30}}
  \put(50,50){\vector(-1,0){17}}

  \qbezier(50,0)(80,20)(80,37)
  \qbezier(80,37)(70,37)(65,50)
  \qbezier(65,50)(60,63)(50,63)
  \put(50,58){\oval(20,10)[tl]} 
  \put(40,58){\vector(0,-1){16}}
  \put(50,42){\oval(20,10)[bl]}
  \qbezier(50,37)(60,37)(65,50)
  \qbezier(65,50)(70,63)(80,63)
  \put(80,58){\oval(20,10)[tr]}
  \put(90,58){\line(0,-1){16}}
  \put(80,42){\oval(20,10)[br]}
  \end{picture}
  }
}

\title{A homological definition of the HOMFLY polynomial}
\author{Stephen Bigelow}
\address{Department of Mathematics,
         University of California at Santa Barbara,
         California 93106, USA}
\email{bigelow@math.ucsb.edu}
\date{September 2006}

\begin{document}

\begin{abstract}
We give a new definition of the knot invariant
associated to the Lie algebra $\su_{N+1}$.
The knot or link must be presented as
the plat closure of a braid.
The invariant is then
a homological intersection pairing between
two submanifolds of a configuration space of points in a disk.
This generalizes previous work on the Jones polynomial,
which is the case $N=1$.
\end{abstract}

\maketitle

\section{Introduction}

The Jones polynomial \cite{vJ85}
was the first of the new generation of knot invariants,
now called ``quantum invariants''.
The two variable HOMFLY polynomial
came soon after \cite{HOMFLY}.
The invariant of type $A_N$ is
a specialization of the HOMFLY polynomial
that is related to the representation theory of $\su_{N+1}$.

Fix an integer $N>1$,
and let $P$ be the invariant of type $A_N$.
This is an invariant of oriented knots and links
that takes values in the $\Zed[q^{\pm 1/2}]$.
It satisfies the following {\em skein relation}.
$$q^{(N+1)/2} P \left( \lminus \right)
  - q^{-(N+1)/2} P \left( \lplus \right)
  = (q^{1/2} - q^{-1/2}) P \left( \lzero \right).$$
Here, the three diagrams represent
three links that are the same except inside a small ball,
where they are as shown.
We can define $P$ to be
the unique invariant that satisfies the above skein relation
and takes the value one for the unknot.

In \cite{sB02},
I presented a definition of the Jones polynomial
as a homological intersection pairing
between a certain pair of manifolds in a configuration space.
The aim of this paper is to give a similar definition
of the invariant of type $A_N$.
The Jones polynomial is the special case $N=1$.
The HOMFLY polynomial can be reconstructed
from the values of all invariants of type $A_N$,
which (perhaps) excuses the title of this paper.

Let $\beta$ be a braid with $2n$ strands.
Orient the strands of $\beta$ in such a way that,
reading from left to right along the bottom of $\beta$,
the orientations are down, up, down, up, and so on.
We require $\beta$ to be such that
reading from left to right along the top of $\beta$,
the orientations are also down, up, down, up, and so on.
Let $\hat{\beta}$ be the plat closure of $\beta$,
obtained by joining adjacent pairs of nodes
at the top and at the bottom of $\beta$.
The orientations on strands of $\beta$
give consistent orientations to the components of $\hat{\beta}$.
Every oriented knot or link
can be obtained as the plat closure of some such braid $\beta$.

The first goal of this paper is to define an invariant $Q(\beta)$.
In Section \ref{sec:configuration},
we define a configuration space $C$.
This is similar to the space $C$ in \cite{sB02}
except that we assign {\em colors}
to the puncture points and the points that make up a configuration.
The colors determine which pairs of points are allowed to coincide,
and how to compute the monodromy of a loop in the configuration space.
In Section \ref{sec:torusball},
we define submanifolds $T$ and $S$ of $C$.
In Section \ref{sec:def},
we define $Q(\beta)$ as an intersection pairing
between $S$ and the image $\beta(T)$ of $T$.

The second goal of this paper is to prove that $Q(\beta) = P(\hat{\beta})$.
In Sections
\ref{sec:heightpreserving},
\ref{sec:bridgepreserving},   and
\ref{sec:markov},
we prove that $Q(\beta)$ is invariant under certain moves.
By a result of Birman \cite{jB76},
this implies that $Q(\beta)$ is
an invariant of the oriented knot or link $\hat{\beta}$.
The more difficult moves require some special tools,
which we develop in Sections \ref{sec:barcodes} and \ref{sec:partialbarcode}.
In Section \ref{sec:skein},
we prove that $Q(\beta)$ satisfies the above skein relation
In Section \ref{sec:conclusion},
we bring these results together to show that $Q(\beta) = P(\hat{\beta})$.

Lawrence gave similar homological definitions
of the Jones polynomial and the invariant of type $A_N$
in \cite{rL93} and \cite{rL96}.
The definition here appears different,
and includes a more precise description
of the relevant manifolds in the configuration space.
Under close examination,
the two approaches might turn out to be the same.

One possible future application of this paper
is to generalize the ideas in \cite{cM04}.
There,
Manolescu gives evidence of a connection between
the definition of the Jones polynomial in \cite{sB02}
and the invariant defined by Seidel and Smith in \cite{SS04}.
Both definitions involve
intersections between submanifolds of configuration spaces.
Seidel and Smith obtain a graded abelian group,
which they conjecture to be
a collapsed version of Khovanov's homology theory.
It would be interesting if the intersection pairing
in \cite{sB02} and in this paper
could be refined to give a graded abelian group.

\rk{Acknowledgements}
This research was partly supported by NSF grant
DMS-0307235 and Sloan Fellowship BR-4124.
I am grateful to Ciprian Manolescu and Dylan Thurston
for their interest and useful conversations.

\section{The configuration space}
\label{sec:configuration}

In this section,
we define the configuration space $C$,
as well as some other terms that will be used throughout the paper.

Let $q$ be a transcendental complex number with unit norm,
and fix a choice of $q^{1/2} \in \Cplx$.
Thus we will work over $\Cplx$ instead of $\Zed[q^{\pm 1/2}]$.
We define the invariant over a more general ring
in Section \ref{sec:conclusion}.

The braid group $B_k$ has many equivalent definitions,
including:
the mapping class group of a $k$-times punctured disk,
the fundamental group of a certain configuration space,
and the group of geometric braids with $k$ strands.
We will move freely between these definitions.
Elements of the mapping class group act on the left,
paths in a configuration space compose from left to right,
and geometric braids read from top to bottom.

Suppose $\pcolor = (c_1,\dots,c_k)$
is a $k$-tuple of elements of $\{0,N+1\}$.
Let $D$ be the unit disk in the complex plane.
Choose points $p_1,\dots,p_k$
ordered from left to right on the real line in $D$.
We will call these {\em puncture points}.
We call $c_i$ the {\em color} of the puncture point $p_i$.
We use the notation $D_\pcolor$ to represent this data.
A braid in $B_k$ induces
a permutation of the puncture points in $D_\pcolor$.
Let the {\em mixed braid group} $B_\pcolor$
be the subgroup of $B_k$
consisting of braids that preserve the colors of the puncture points.

Suppose $\mcolor = (c'_1,\dots,c'_m)$
is an $m$-tuple of elements of $\{1,\dots,N\}$.
We now define the configuration space $C_{\mcolor}(D_{\pcolor})$.
First, let $\tilde{C}$ be the set of all
$m$-tuples $(x_1,x_2,\dots,x_m)$ of points in $D$ such that
\begin{itemize}
\item if $1 \le i < j \le m$ and $|c'_i - c'_j| \le 1$
      then $x_i \neq x_j$, and
\item if $1 \le i \le m$, $1 \le j \le k$, and $|c'_i - c_j| = 1$
      then $x_i \neq p_j$.
\end{itemize}
Now let $W$ be the group of
permutations of $\{1,\dots,m\}$
such that $c'_i = c'_{w(i)}$ for all $i=1,\dots,m$.
Let $C_\mcolor(D_\pcolor)$
be the quotient of $\tilde{C}$ by the induced action of $W$.

Thus a point in $C_\mcolor(D_\pcolor)$
is a configuration of $m$ points in $D$,
which we call {\em mobile points}.
These mobile points have colors given by $\mcolor$.
Two mobile points of the same color are indistinguishable.
A mobile point may coincide with
a puncture point or another mobile point
if and only if their colors differ by at least two.

We will represent elements of $\pi_1(C)$ using braids as follows.
Let $\pcolor + \mcolor$ denote the concatenation
$$\pcolor + \mcolor = (c_1,\dots,c_k,c'_1,\dots,c'_m).$$
Let $G$ be group of those mixed braids in $B_{\pcolor + \mcolor}$
whose first $k$ strands are straight.
Then $\pi_1(C)$ is the quotient of $G$ obtained by equating
any two braids that differ by a sequence of crossing changes
involving pairs of strands whose colors differ by at least two.
Thus we can represent an element of $\pi_1(C)$ by a braid in $G$.
We will put the straight strands corresponding to puncture points
in whatever position is convenient,
and not necessarily on the left.

Let $\pcolor$ be the $2n$-tuple
$$\pcolor = (0,N+1,0,N+1,\dots,0,N+1).$$
Note that our braid $\beta$ is an element of $B_\pcolor$,
where the strands of with color $N+1$ are oriented upwards,
and strands with color $0$ are oriented downwards.
Let $m = Nn$,
and let $\mcolor$ be the $m$-tuple
$$\mcolor = (1,2,\dots,N,1,2,\dots,N,\dots,1,2,\dots,N).$$
Let $C$ denote the configurations space $C_\mcolor(D_\pcolor)$.

We now define a homomorphism
$$\rho_\mcolor \co \pi_1(C) \to \{\pm q^k \suchthat k \in \Zed\}.$$
Suppose $g$ is an element of $\pi_1(C)$.
Represent $g$ by a braid diagram.
To every positive crossing in this braid diagram,
associate the term
\begin{itemize}
\item $-q^{-1}$ if it involves two strands of the same color,
\item $q^{1/2}$ if it involves two strands whose colors differ by one,
\item $1$ otherwise.
\end{itemize}
To every negative crossing,
associate the reciprocal of
the term associated to the analogous positive crossing.
Let $\rho_\mcolor(g)$ be the product of the terms
associated to the crossings of the braid diagram.
Note that the exponent of $q$ in $\rho_\mcolor(g)$ is an integer,
since there must be an even number of crossings
involving strands whose colors differ by one.

Next we define a homomorphism
$$\rho_\pcolor \co B_\pcolor \to \{\pm q^{k/2} \suchthat k \in \Zed\}.$$
Suppose $g$ is an element of $B_\pcolor$.
Represent $g$ by a braid diagram.
To every positive crossing in $g$,
associate the monomial
\begin{itemize}
\item $q^{N/2}$ if it involves two strands of the same color,
\item $q^{-(N+1)/2}$ if it involves two strands of different colors.
\end{itemize}
To every negative crossing,
associate the reciprocal of
the term associated to the analogous positive crossing.
Let $\rho_\pcolor(g)$ be the product of the monomials
associated to the crossings of $g$.

\section{A torus and a ball}
\label{sec:torusball}

The aim of this section is to define
an immersion $\Phi$ from an $m$-dimensional torus to $C$,
and an embedding $\Psi$ from an open $m$-ball to $C$.
Until otherwise stated,
we assume that $n=1$,
and hence that $\pcolor = (0,N+1)$ and $\mcolor = (1,2,\dots,N)$.

Let $S^1$ be the unit circle centered at the origin in the complex plane,
and let $T$ be the product of $N$ copies of $S^1$.
Let $A$ and $B$ be the intersections of $S^1$
with the closed upper and lower half planes respectively.

\begin{figure}
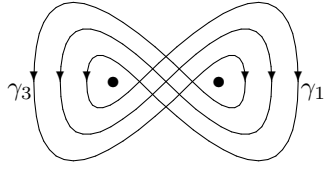

$$\eights$$
\caption{Figures of eight in the case $N=3$.}
\label{fig:eights}
\end{figure}

Let $\gamma_1,\dots,\gamma_N \co S^1 \to D$
be figures of eight as shown in Figure \ref{fig:eights}.
Assume $\gamma_i$ is parametrized so that
$\gamma_i|_A$ is a loop that winds counterclockwise around $p_1$,
and $\gamma_i|_B$ is a loop that winds clockwise around $p_2$.
Thus the loops $\gamma_i(A)$ are concentric loops around $p_1$,
and the loops $\gamma_i(B)$ are concentric loops around $p_2$.
We assume that the points $\gamma_i(1)$ are all on the real line,
and
$$p_1 < \gamma_1(1) < \dots < \gamma_N(1) < p_2.$$

\subsection{The case $N=2$}

We now define $\Phi \co T \to C$ in the case $N=2$.
The most difficult part of $\Phi$ is given by the following lemma.

\begin{lem}
\label{lem:tricky}
There is an immersion $\Phi_1 \co B \times A \to C$ such that
    $$\Phi_1 |_{\partial(B \times A)} =
      (\gamma_1 \times \gamma_2) |_{\partial(B \times A)}.$$
\end{lem}

\begin{proof}
It suffices to show that the loop
$(\gamma_1 \times \gamma_2) |_{\partial(B \times A)}$
is null-homotopic in $C$.
This loop is the commutator of the loops $\alpha$ and $\beta$,
where $\alpha \co A \to C$ is given by
$$\alpha(s) = (\gamma_1(1), \gamma_2(s)),$$
and $\beta \co B \to C$ is given by
$$\beta(s) = (\gamma_1(s), \gamma_2(1)).$$
We can represent $\alpha$ by a braid as follows.
$$\alpha = \baab$$
Here,
the strands are colored $0$, $1$, $2$, and $3$, from left to right.
The strands on the far left and the far right represent the puncture points.

Recall that
the strand of color $2$ may pass through the strand of color $0$.
Thus $\alpha$ is homotopic relative to endpoints to the loop $\alpha'$
represented by the braid as follows.
$$\alpha' = \bb$$
Similarly,
$\beta$ is homotopic relative to endpoints to $\beta' = (\alpha')^{-1}$.
Then $\alpha'$ and $\beta'$ obviously commute,
thus completing the proof.
\end{proof}

We can now define $\Phi \co T \to C$ as follows.
$$\Phi(s_1,s_2) =
\left\{
\begin{array}{ll}
\Phi_1(s_1,s_2)               & \mbox{if $(s_1,s_2) \in B \times A$,} \\
(\gamma_1(s_1),\gamma_2(s_2)) & \mbox{otherwise.}
\end{array}
\right.
$$

This completes the definition of $\Phi$ when $n=1$ and $N=2$.
We can choose $\Phi_1$ to have some properties
that will be useful later.

\begin{lem}
\label{lem:trickytricky}
The function $\Phi_1$ in the previous lemma
can be chosen so that for every $(x_1,x_2)$ in its image,
\begin{itemize}
\item $x_1$ lies in the closed disk bounded by $\gamma_1(B)$,
\item $x_2$ lies in the closed disk bounded by $\gamma_2(A)$, and
\item at least one of $x_1$ and $x_2$ lies in
      the intersection of these two disks.
\end{itemize}
\end{lem}

\begin{proof}
Let $C'$ be the set of points $(x_1,x_2) \in C$
satisfying the three requirements of the lemma.
Let $C''$ be the set of points $(x_1,x_2) \in C$ such that
$x_1$ and $x_2$ both lie in the intersection of the closed disks
bounded by $\gamma_1(B)$ and $\gamma_2(A)$.
Let $\alpha$, $\alpha'$, $\beta$ and $\beta'$
be as in the proof of the previous lemma.

Any reasonable choice of homotopy
from $\alpha$ to $\alpha'$ relative to endpoints
will lie in $C'$.
Further,
we can assume that $\alpha'$ lies in $C''$.
Similarly,
we can assume that
the homotopy from $\beta$ to $\beta'$ lies in $C'$,
and $\beta'$ lies in $C''$.
The commutator of $\alpha'$ and $\beta'$ is null homotopic
as a loop in $C''$.
\end{proof}

\subsection{General values of $N$}

We now define $\Phi \co T \to C$ for general values of $N$.
We will use functions
$$\Phi_1,\dots,\Phi_{N-1} \co B \times A \to D \times D,$$
similar to $\Phi_1$ for the case $N=2$.
Specifically,
$$\Phi_i |_{\partial (B \times A)}
  = (\gamma_i \times \gamma_{i+1}) |_{\partial (B \times A)},$$
and for all $(x_i,x_{i+1})$ in the image of $\Phi_i$,
\begin{itemize}
\item $x_i \neq x_{i+1}$,
\item $x_i$ lies in the closed disk bounded by $\gamma_i(B)$,
\item $x_{i+1}$ lies in the closed disk bounded by $\gamma_{i+1}(A)$, and
\item at least one of $x_i$ and $x_{i+1}$ lies in
      the intersection of these two disks.
\end{itemize}

Suppose $(s_1,\dots,s_N) \in T$.
For $i = 1,\dots,N$, let $x_i$ be as follows.
\begin{itemize}
\item if $s_i,s_{i-1} \in A$ then $x_i = \gamma_i(s_i)$,
\item if $s_i \in A$ and $s_{i-1} \in B$
      then $x_i$ is the second coordinate of $\Phi_{i-1}(s_{i-1},s_i)$,
\item if $s_i,s_{i+1} \in B$ then $x_i = \gamma_i(s_i)$, and
\item if $s_i \in B$ and $s_{i+1} \in A$
      then $x_i$ is the first coordinate of $\Phi_i(s_i,s_{i+1})$.
\end{itemize}
Here,
for convenience,
we take $s_0$ to be a point in $A \setminus B$
and $s_{N+1}$ to be a point in $B \setminus A$.
Let $\Phi(s_1,\dots,s_N) = (x_1,\dots,x_N)$.

We must show that $\Phi$ is a well defined map from $T$ to $C$.
First note that
if two or more of the conditions apply in the definition of $x_i$
then they all give the value $x_i = \gamma_i(s_i)$.
Next note that $x_1 \neq p_1$,
since either $x_1 = \gamma_1(s_1)$
or $x_1$ is the first coordinate of $\Phi_1(s_1,s_2)$.
Similarly,
$x_N \neq p_2$.
It remains to show that $x_i \neq x_{i+1}$ for all $i=1,\dots,N-1$.
There are several cases to check.

First, suppose $s_i \in A$.
Then either $x_i = \gamma_i(s_i)$
or $x_i$ is the second coordinate of $\Phi_{i-1}(s_{i-1},s_i)$.
Also, either $x_{i+1} = \gamma_{i+1}(s_{i+1})$
or $x_{i+1}$ is the first coordinate of $\Phi_{i+1}(s_{i+1},s_{i+2})$.
In all cases,
$x_i$ lies in the disk bounded by $\gamma_i(A)$,
and $x_{i+1}$ does not.
Thus $x_i \neq x_{i+1}$.

The case $s_{i+1} \in B$ is similar.

Finally, if $s_i \in B$ and $s_{i+1} \in A$
then $(x_i,x_{i+1}) = \Phi_i(s_i,s_{i+1})$,
so $x_i \neq x_{i+1}$.

This completes the proof that $\Phi$ is a well defined map from $T$ to $C$.
It also has the following important property.

\begin{lem}
\label{lem:toruslifts}
$\rho_\mcolor \circ \Phi_*(\pi_1(T)) = \{1\}$.
\end{lem}
\begin{figure}
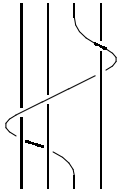

$$\babCC$$
\caption{A braid representing $g_2$ when $N=2$.}
\label{fig:eightbraid}
\end{figure}
\begin{proof}
For $i=1,\dots,N$,
let $g_i \co S^1 \to T$ be the map
$$g_i(s) = (\gamma_1(1), \dots, \gamma_{i-1}(1), \gamma_i(s),
            \gamma_{i+1}(1), \dots, \gamma_N(1)).$$
These loops generate $\pi_1(T)$.
We must show that $\rho_\mcolor(g_i) = 1$.
Represent $g_i$ by a mixed braid with $N+2$ strands.
Every strand is straight except for the strand with color $i$,
which describes a figure of eight.
See Figure \ref{fig:eightbraid}.
There are two positive crossings
that involve a pair of strands with colors $i$ and $i-1$,
and two negative crossings
that involve a pair of strands with colors $i$ and $i+1$.
Thus $\rho_\mcolor(g_i) = 1$.
\end{proof}

\subsection{Working with the immersed torus}
\label{subsec:X}

We now describe how to partition $\Phi(T)$ into two parts,
one of which is easy to work with,
and the other of which can be safely ignored.

Let $X$ be the the intersection of the disks
bounded by $\gamma_1(B)$ and $\gamma_N(A)$.
Let $C_X$ be the set of points in $C$
that include a mobile point in $X$.
Let $T_X$ be the set of $(s_1,\dots,s_N) \in T$
such that $s_i \in B$ and $s_{i+1} \in A$
for some $i=1,\dots,N-1$.
Then $\Phi(T_X)$ lies in $C_X$.
On the other hand,
$\Phi(T \setminus T_X)$ is a disjoint union of $N+1$ embedded $N$-balls.
In practice,
we can usually take $X$ to be small,
ignore $\Phi(T_X)$,
and restrict our attention to $\Phi(T \setminus X)$.

From now on we will omit any reference to $\Phi$,
and simply treat $T$ as
oriented $N$-dimensional submanifold of $C$.

\subsection{A basepoint}

Choose points $t_1,\dots,t_N$ in the disk such that
\begin{itemize}
\item $t_i \in \gamma_i(B)$,
\item $t_i$ is below the real line,
\item the real parts of $t_1,\dots,t_N$
      are in increasing order
      and lie between $\gamma_N(1)$ and $p_2$.
\end{itemize}
Let $\mathbf{t} = (t_1,\dots,t_N)$.
This will be our basepoint of $T$.

For $i=1,\dots,N$,
let $\tau_i \co I \to D$ be a vertical edge
from a point on the lower half of $\partial D$ up to $t_i$.
Let $\tau \co I \to C$ be the path
$$\tau(s) = (\tau_1(s),\dots,\tau_N(s)).$$
Let $\mathbf{x} = \tau(0)$.
This will be our basepoint for $C$.
Thus $\tau$ is a path from the basepoint $\mathbf{x}$ of $C$
to the basepoint $\mathbf{t}$ of $T$.

\subsection{A ball}

Let
$$S = \{(s_1,\dots,s_N) \in \Real^N \suchthat 0 < s_1 < \dots < s_N < 1\}.$$
This is an open $N$-ball.
Let $\gamma \co I \to D$ be the straight edge from $p_1$ to $p_2$.
Let $\Psi \co S \to C$ be the embedding
$$\Psi(s_1,\dots,s_N) = (\gamma(s_1), \dots, \gamma(s_N)).$$
For $i=1,\dots,N$,
let $\zeta_i \co I \to D$ be a vertical edge
from $x_i$ to a point on $\gamma$.
Let $\zeta \co I \to C$ be the map
$$\zeta(s) = (\zeta_1(s),\dots,\zeta_N(s)).$$
Let $\mathbf{s} = \zeta(1)$.
This will be our basepoint for $S$.
Thus $\zeta$ is an path from the basepoint $\mathbf{x}$ of $C$
to the basepoint $\mathbf{s}$ of $S$.

From now on we will omit any reference to $\Psi$,
and simply treat $S$ as
an oriented $N$-dimensional submanifold of $C$.

\subsection{General values of $n$}

We now define $T$ and $S$ for general values of $n$.

Let $C_1 = C_{(1,\dots,N)}(D_{(0,N+1)})$.
This is the configuration space in the case $n=1$.
Note that $D_\pcolor$ can be obtained
by gluing together $n$ copies of $D_{(0,N+1)}$ side by side.
This defines an embedding
from the product of $n$ copies of $C_1$ into $C$.

Let $T$ be the product of $n$ copies of the immersed $N$-torus in $C_1$
as defined in the case $n=1$.
Also let $\tau$ be the product of
$n$ copies of the path in $C_1$.
This is a path from a basepoint $\mathbf{x}$ of $C$
to a basepoint $\mathbf{t}$ of $T$.

Define an open $m$-ball $S$
and a path $\zeta$ from $\mathbf{x}$ to a basepoint $\mathbf{s}$ of $S$
similarly,
by taking a product of $n$ copies of the versions when $n=1$.

\section{Definition of the invariant}
\label{sec:def}

The aim of this section is to define the invariant $Q(\beta)$.
We give two equivalent definitions
of an intersection pairing $\langle S,\beta(T) \rangle$.
The first gives an explicit method of computation,
and the second uses a more abstract homological approach.
We then define $Q(\beta)$
to be a renormalization of $\langle S,\beta(T) \rangle$.

\subsection{An intersection pairing}

We can represent $\beta$ by a homeomorphism from $D$ to itself
that preserves the colors of the puncture points.
This induces a homeomorphism from $C$ to itself,
which we also call $\beta$.

Note that $S$ and $\beta(T)$
are immersed $m$-manifolds in the $(2m)$-manifold $C$.
By applying a small isotopy we can assume that 
they intersect transversely at a finite number of points.
For each such intersection point $\mathbf{y}$,
let $\epsilon_\mathbf{y}$ be the sign
of the intersection at $\mathbf{y}$,
and let $\xi_\mathbf{y}$ be
the composition of the following paths in order.
\begin{itemize}
\item $\beta \circ \tau$,
\item an path in $\beta(T)$ from $\beta(\mathbf{t})$ to $\mathbf{y}$,
\item an path in $S$ from $\mathbf{y}$ to $\mathbf{s}$,
\item $\overline{\zeta}$.
\end{itemize}
Let
$$\langle S,\beta(T) \rangle
  = \sum \epsilon_\mathbf{y} \rho_\mcolor(\xi_\mathbf{y}),$$
where the sum is taken over all $\mathbf{y} \in S \cap \beta(T)$.

We now describe how one could
use this definition to explicitly compute
$\langle S,\beta(T) \rangle$ for a given $\beta$.
The computation is complicated,
and impractical in all but the simplest examples.
However it might provide an aid to understanding,
and some aspects of it will be used later in the paper.

Recall that $T$ was defined to be the product of
$n$ copies of an $N$-dimensional torus.
Call these tori $T_1,\dots,T_n$.
Corresponding to each torus
is a small disk $X$ as defined in Section \ref{subsec:X}.
We assume that the images of these disks under $\beta$
are disjoint from the intervals $[p_{2j-1},p_{2j}]$ used to define $S$.

We first describe how to recognize a point $\mathbf{y}$
in the intersection of $S$ and $\beta(T)$.
Note that $\mathbf{y}$ lies in $S$ if and only if
every interval $[p_{2i-1},p_{2i}]$
contains $N$ of the mobile points of $\mathbf{y}$,
having colors $1,\dots,N$,
reading from left to right.

Now $\mathbf{y}$ lies on $\beta(T)$
if and only if the following conditions hold for each $i=1,\dots,N$.
Let $\gamma_1,\dots,\gamma_N$ be the figures of eight
used to define $T_i$.
Then,
for every $j=1,\dots,N$,
$\mathbf{y}$ must include
one mobile point of color $j$ on $\beta(\gamma_j)$.
This lies in one of the two loops that make up $\beta(\gamma_j)$.
Taking the corresponding loops for all $j$,
we must have the innermost $N_1$ loops around one of the puncture points,
and the innermost $N_2$ loops around the other,
for some $N_1$ and $N_2$ with $N_1+N_2=N$.

Next we compute a braid diagram representing the path $\xi_\mathbf{y}$.
We do this first in the case $n=1$.
The two strands corresponding to puncture points will always be straight.
For $i=1,\dots,N$,
the strand of color $i$ describes a path
along $\beta(\tau_i)$ and then along $\beta(\gamma_i)$
to the mobile point that lies on this figure of eight.
The order these strands follow these paths
is not important
except that the paths along $\beta(\gamma_i)$
must be performed in order $i=1,\dots,N$.
Note that the last half of $\xi_\mathbf{y}$,
which lies in $S$ and $\zeta$,
contributes no crossings to the braid $\xi_\mathbf{y}$.

The case $n > 1$ is basically the same.
The $2n$ strands corresponding to puncture points are straight.
Each of the remaining $m$ strands describes a path
along the copies of $\beta(\tau_i)$ and $\beta(\gamma_i)$
corresponding to the appropriate torus $T_j$.
The order is not important
except that within each torus $T_j$,
the paths along $\beta(\gamma_i)$
must be performed in order $i=1,\dots,N$.

We now compute the sign $\epsilon_\mathbf{y}$.
Each mobile point of $\mathbf{y}$ is a point of intersection
between some edge from $p_{2j-1}$ to $p_{2j}$
and the image under $\beta$
of one of the figures of eight used to define $T$.
Determine the sign of this intersection,
taking the oriented edge first,
and the oriented figure of eight second.
Then $\epsilon_\mathbf{y}$ is the product of
the signs of the intersections at the mobile points of $\mathbf{y}$,
multiplied by the sign of the permutation of the mobile points
induced by the loop $\xi_\mathbf{y}$.

This completes the computation of $\langle S,\beta(T) \rangle$.
By Lemma \ref{lem:toruslifts},
$\rho_\mcolor(\xi_\mathbf{y})$
does not depend on the choice of path in $\beta(T)$.
It remains to check that the sum
is invariant under isotopy of $\beta$.
One could do this by checking invariance under certain moves.
However the real reason $\langle S,\beta(T) \rangle$ is well defined
is that it computes the homological intersection pairing described below.

\subsection{A homological definition}
\label{subsec:homdef}

We now define some homology modules of $C$.

Let $\Local$ be the flat
complex line bundle over $C$
with monodromy given by $\rho_\mcolor$.
Let $\langle \cdot,\cdot \rangle$
be the sesquilinear inner product on $\Cplx$
given by $\langle x,y \rangle = \bar{x}y$.
This inner product is preserved by the monodromy of $\Local$,
so it gives a well-defined inner product on the fiber of $\Local$ at any point.
In other words,
$\Local$ is a Hilbert line bundle.
Topologists may prefer to give $\Cplx$ the discrete topology
and think of $\Local$ as a covering space of $C$.
Each fiber of this covering space has the structure of
a $1$-dimensional Hilbert space,
and these structures are locally consistent.

Let $H_m(C;\Local)$ denote
the $m$-dimensional homology of $C$ with local coefficients.
For a definition of homology with local coefficients,
see, for example, \cite[Section 3H]{aH02}.
The idea is the same as singular homology with module coefficients,
except that the coefficient of a simplex
is a lift of that simplex to $\Local$.

Let $H_m^\lf(C;\Local)$ denote the $m$-dimensional
{\em locally finite} homology of $C$ with local coefficients
(also called {\em Borel-Moore} homology).
For a definition of locally finite homology,
see, for example, \cite[Exercise 3H.6]{aH02}.
Briefly,
the idea is to allow infinite sums of simplices with local coefficients,
as long as every compact set in $C$ meets only finitely many simplices.

From now on,
all homology modules will be assumed to use coefficients in $\Local$.
For example, we will write $H_m(C)$ to mean $H_m(C;\Local)$.
We also use relative versions of these homology theories.
Recall the following basic theorems.

\begin{thm}[Poincar\'e-Lefschetz Duality]
$H^m(C)$ and $H_m^\lf(C,\partial C)$ are isomorphic.
\end{thm}

\begin{thm}[The Universal Coefficient Theorem]
$H^m(C)$ and $\Hom(H_m(C),\Cplx)$ are conjugate-isomorphic.
\end{thm}

These theorems imply that
$H_m^\lf(C,\partial C)$ and $\Hom(H_m(C),\Cplx)$ are conjugate-isomorphic.
Thus there is a sesquilinear pairing
$$\langle \cdot,\cdot \rangle  \co
  H_m^\lf(C,\partial C) \times H_m(C) \to \Cplx.$$
The precise definition of this pairing
follows from the more explicit statements of
Poincar\'e-Lefschetz duality and the universal coefficient theorem,
which give the definitions of the isomorphisms.

Let $\beta_*$ be the automorphism of $\pi_1(C)$ induced by $\beta$.
It is not too hard to show that $\rho_\mcolor \circ \beta_* = \rho_\mcolor$.
Thus $\beta$ lifts to an action on $\Local$.
Choose this lift to act as the identity
on the fiber over the basepoint $\mathbf{x}$.
Thus there are induced actions of $\beta$
on $H_m(C)$, $H_m(C,\partial C)$ and $H_m^\lf(C)$.
By abuse of notation,
we use $\beta$ to denote every one of these induced actions.

For the rest of this paper,
fix an identification of the fiber over $\mathbf{x}$ with $\Cplx$.
Let $\tilde{\tau}$ be the lift of $\tau$ to $\Local$
starting at the element $1$ of the fiber over $\mathbf{x}$.
By Lemma \ref{lem:toruslifts},
we can lift $T$ to an immersed torus $\tilde{T}$ in $\Local$
such that $\tilde{T}$ contains $\tilde{\tau}(1)$.
This determines an element of $H_m(C)$,
which we also denote by $T$.

Similarly,
$\zeta$ determines a lift $\tilde{S}$ of $S$ to $\Local$.
Let $S$ denote
the open $m$-ball,
the corresponding element of $H_m^\lf(C)$,
and also the corresponding element of $H_m^\lf(C,\partial C)$.
Then $\langle S,\beta(T) \rangle$
is the sesquilinear pairing of
$S \in H_m^\lf(C,\partial C)$ and $\beta(T) \in H_m(C)$.

We list some properties of the pairing $\langle \cdot,\cdot \rangle$.
\begin{itemize}
\item It is the same as the previous more computational definition,
\item it is sesquilinear
      (conjugate-linear in the first entry and linear in the second),
\item it is invariant under the action of $B_\pcolor$,
\item it has the following symmetry property:
      if $v_1,v_2 \in H_m(C)$
      and $v'_1,v'_2$ are their images in $H_m^\lf(C,\partial C)$ then
      $\langle v'_1,v_2 \rangle = (-1)^m \overline{\langle v'_2,v_1 \rangle}$.
\end{itemize}
These all follow from standard homology theory.

As an aside, note that
it might be possible to obtain a unitary representation of $B_\pcolor$
with some more work along these lines.
Compare the result of Budney \cite{rB05}
that the Lawrence-Krammer representation
is negative-definite Hermitian.

\subsection{Definition of the invariant}

We are finally ready to define the invariant $Q(\beta)$.
Let
$$[N+1] = \frac{q^{(N+1)/2} - q^{-(N+1)/2}}{q^{1/2} - q^{-1/2}}.$$
This is the {\em quantum integer} corresponding to $N+1$.
Then let
$$Q(\beta) = \frac{\rho_\pcolor(\beta)}{[N+1]q^{m/2}}
             \langle S,\beta(T) \rangle.$$
The main result of this paper is that $Q(\beta) = P(\hat{\beta})$.

\section{Height-preserving isotopy}
\label{sec:heightpreserving}

For all $i=1,\dots,b-1$,
let $\sigma'_i = \sigma_{2i} \sigma_{2i+1} \sigma_{2i-1} \sigma_{2i}$.
The aim of this section is to prove the following.

\begin{lem}
\label{lem:heightpreserving}
$Q(\sigma_1^2 \beta) = Q(\beta \sigma_1^2)
 = Q(\sigma'_i \beta) = Q(\beta \sigma'_i) = Q(\beta)$.
\end{lem}

Assume the plat closure $\hat{\beta}$
is defined so that all maxima are at the same height
and all minima are at the same height.
Then the above lemma is equivalent to the statement that
$Q(\beta)$ is invariant under height preserving isotopy of $\hat{\beta}$.
We will not use this formulation,
but mention it by way of motivation.

\begin{clm}
$Q(\sigma_1^2 \beta) = Q(\beta)$.
\end{clm}

\begin{proof}
We have
$$\rho_\pcolor(\sigma_1^2 \beta) = q^{-(N+1)} \rho_\pcolor(\beta).$$
By this and the properties of the sesquilinear pairing,
it suffices to show
$$\sigma_1^2 S = q^{N+1}S.$$
We can choose the function $\sigma_1^2$ to act as the identity on
the subset $S$ of $C$.
It remains to show that $\sigma_1^2$ acts as multiplication by $q^{N+1}$
on the fiber over $\mathbf{s}$.

\begin{figure}
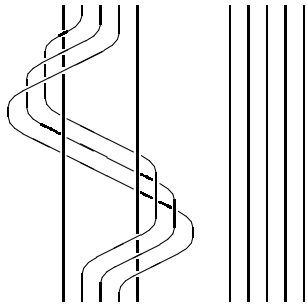

$$\fulltwist$$
\caption{$(\sigma_1^2 \zeta)\cdot(\overline{\zeta})$ when $n=2$ and $N=3$.}
\label{fig:fulltwist}
\end{figure}

Let $\xi$ be the concatenation of the paths
$\sigma_1^2 \zeta$ and $\overline{\zeta}$.
This is represented by a braid in which
strands of colors $1,\dots,N$
make a positive full twist around two with colors $0$ and $N+1$.
Figure \ref{fig:fulltwist} shows this braid when $n=1$ and $N=3$.
Then
$$\rho_\mcolor(\xi) = (q^{1/2})^{2N+2} = q^{N+1}.$$
Thus $\sigma_1^2(S) = q^{N+1} S$,
as required.
\end{proof}

\begin{clm}
$Q(\sigma'_i \beta) = Q(\beta)$.
\end{clm}

\begin{proof}
We have
$$\rho_\pcolor(\sigma'_i \beta) = q^{-1} \rho_\pcolor(\beta).$$
By this and the properties of the sesquilinear pairing,
it suffices to show
$$\sigma'_i S = qS.$$
We can choose the function $\sigma'_i$ to act as the identity on
the subset $S$ of $C$.
It remains to show that $\sigma'_i$ acts as multiplication by $q$
on the fiber over $\mathbf{s}$.

\begin{figure}
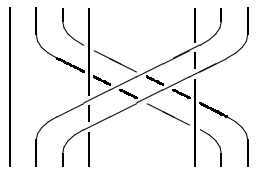

$$\bigx$$
\caption{$(\sigma'_1 \zeta)\cdot(\overline{\zeta})$ when $n=2$ and $N=2$.}
\label{fig:bigx}
\end{figure}

Let $\xi_i$ be the concatenation of the paths
$\sigma'_i \zeta$ and $\overline{\zeta}$.
This is represented by a braid in which
two collections of $N$ parallel strands of colors $1,\dots,N$
form a large letter $X$
enclosing two strands of colors $0$ and $N+1$.
Figure \ref{fig:bigx} shows this braid when $N=2$, $i=1$ and $n=2$.
Then
$$\rho_\mcolor(\xi_i) = (q^{-1})^N (q^{1/2})^{2N+2} = q.$$
Thus $\sigma'_i(S) = q S$,
as required.
\end{proof}

It remains to show that
$Q(\beta \sigma_1^2) = Q(\beta)$ and $Q(\beta \sigma'_i) = Q(\beta)$.
It suffices to show that $\sigma_1^2 T = q^{N+1}T$ and $\sigma'_i(T) = q T$.
The proof of these identities
is the same as the proof of the analogous identities for $U$
given in the previous two claims.
This completes the proof of the lemma. 

\section{Barcodes}
\label{sec:barcodes}

Before we prove the invariance of $Q(\beta)$ under other moves,
we will look more closely at $H_m^\lf(C)$ and $H_m(C,\partial C)$.
In the process,
we will introduce a useful tool I call a {\em barcode}.

\subsection{A basis for $H_m^\lf(C)$}

Let $C_\Real$ be the set of points in $C$
that are configurations of points on the real line in $D$.

\begin{lem}
\label{lem:getreal}
The map $H_m^\lf(C_\Real) \to H_m^\lf(C)$ induced by inclusion
is an isomorphism.
\end{lem}

For the details of the proof,
see \cite[Lemma 3.1]{sB04}.
The idea is to vertically ``squash''
configurations of points in the disk
to configurations of points in the real line.
The only difficulty is that
a configuration may contain two mobile points,
or a mobile point and a puncture point,
that are mapped to the same point on the real line,
although their colors differ by at most one.
Such a configuration would be ``sent to infinity''
as it is squashed to the real line.
Since we are using locally finite homology,
this does not pose a serious problem.

We now enumerate the components of $C_\Real$.

\begin{defn}
A {\em code sequence} is a permutation of the sequence $\pcolor + \mcolor$
that contains $\pcolor$ as a subsequence.
\end{defn}

Suppose $S'$ is a connected component of $C_\Real$.
Choose a point $\mathbf{y} = (y_1,\dots,y_m)$ in $S'$
such that $y_1,\dots,y_m$ are
distinct from each other and from the puncture points.
Let $\mathbf{c} = (c_1,\dots,c_{m+n})$
be the sequence of colors of mobile points and puncture points,
reading from left to right on the real line.
Then $\mathbf{c}$ is a code sequence.
We say $\mathbf{c}$ {\em represents} $S'$.

Suppose $i$ is such that
at least one of $c_i$ and $c_{i+1}$ is in $\{1,\dots,N\}$
and $|c_i - c_{i+1}| \ge 2$.
Then we can exchange $c_i$ and $c_{i+1}$ in $\mathbf{c}$
without altering the connected component of $C_\Real$ it represents.
This corresponds to
moving a mobile point
through another mobile point or a puncture,
provided their colors permit this.
We say two code sequences are {\em equivalent}
if they are related by a sequence of such transpositions.
The equivalence classes of code sequences
enumerate the connected components of $C_\Real$.

\begin{defn}
A code sequence is {\em trivial}
if it is equivalent to a code sequence
whose first or last entry lies in $\{1,\dots,N\}$.
\end{defn}

Suppose $S'$ is the connected component of $C_\Real$
corresponding to a code sequence $\mathbf{c}$.
If $\mathbf{c}$ is trivial
then $S'$ contains a point $(y_1,\dots,y_m)$
such that $y_1$ or $y_m$ lies on $\partial D$.
In this case, $S'$ is homeomorphic to the upper half space in $\Real^m$,
so $H_m^\lf(S') = 0$.
If $\mathbf{c}$ is not trivial
then every point in $S'$ is
a configuration of points between $p_1$ and $p_{2n}$.
In this case,
$S'$ is homeomorphic to an open $m$-ball,
so $H_m^\lf(S') = \Cplx$.

For every nontrivial code sequence $\mathbf{c}$,
choose a nonzero element of $H_m^\lf(S')$,
where $S'$ is the corresponding component of $C_\Real$.
By Lemma \ref{lem:getreal},
this gives a basis for $H_m^\lf(C)$.
To define this basis precisely,
we would need to specify
an orientation and a lift to $\Local$
for every component of $C_\Real$.
In practice,
it often suffices to specify an element of $H_m^\lf(C)$
up to multiplication by a nonzero scalar.

\subsection{A basis for $H_m(C,\partial C)$}

Let
$$\langle \cdot,\cdot \rangle' \co
  H_m^\lf(C) \times H_m(C,\partial C) \to \Cplx$$
be the nondegenerate sesquilinear pairing defined
using the more general version of the Poincar\'e-Lefschetz Duality.
We define a basis of $H_m(C,\partial C)$
that is dual to our basis of $H_m^\lf(C)$ with respect to this pairing.

Let $E_1,\dots,E_m$ be properly embedded vertical edges in $D$
that are disjoint from each other and from the puncture points.
The product of these edges
is a properly embedded closed $m$-ball $Z$ in $C$.
Let $\mathbf{c} = (c_1,\dots,c_{m+n})$
be the sequence of colors of vertical edges or puncture points,
reading from left to right.
This is a code sequence.

Any nonzero lift of $Z$ to $\Local$
represents an element of $H_m(C,\partial C)$.
By abuse of notation,
we will use $Z$ to denote both the embedded $m$-ball
and a corresponding element of $H_m(C,\partial C)$,
and call either of these
the {\em barcode} corresponding to the code sequence $\mathbf{c}$.

Two equivalent code sequences
will give rise to the same barcode in $H_m(C,\partial C)$,
up to the choices of lifts to $\Local$.
If $\mathbf{c}$ is trivial
then any barcode corresponding to $\mathbf{c}$ is zero.
Choose a nonzero barcode corresponding to
each nontrivial code sequence $\mathbf{c}$.
I claim that these form a basis for $H_m(C,\partial C)$.

Suppose $S'$ is a component of $C_\Real$
and $Z$ is a barcode.
If $S'$ and $Z$ correspond to
the same nontrivial code sequence
then they intersect at one point,
so we can choose our lifts and orientations so that
$\langle S',Z \rangle' = 1$.
On the other hand,
if $S'$ and $Z$ correspond to different nontrivial code sequences
then they do not intersect,
so $\langle S'Z \rangle' = 0$.
Thus we have a basis of $H_m(C,\partial C)$
that is dual to our basis for $H_m^\lf(C)$.

\subsection{Images of $T$}
\label{subsec:images}

Using the above bases,
we now compute the image of $T$ in $H_m^\lf(C)$,
and also in $H_m(C,\partial C)$ in the case $n=1$.

If $n=1$ then the unique nontrivial code sequence is $(0,1,\dots,N+1)$.
Let $Z$ be the corresponding barcode.
We specify an orientation and lift of $Z$ as follows.
Take $Z$ to be the product of edges of colors $1,\dots,N$ in order,
with each edge oriented upwards.
We can assume that $Z$ contains the basepoint $\mathbf{x}$.
Choose the lift of $Z$ to $\Local$
that contains the point $1$ in the fiber over $\mathbf{x}$.
Note that $\langle S,Z \rangle' = 1$.

\begin{lem}
\label{lem:TS}
The image of $T$ in $H_m^\lf(C)$ is $(q-1)^m S$.
\end{lem}

\begin{proof}
By the construction of $T$,
it suffices to prove this lemma in the case $n=1$.

There is only one nontrivial code sequence,
so $T$ is some scalar multiple of $S$ in $H_m^\lf(C)$.
It remains to show that
$$\langle T,Z \rangle' = \overline{(q-1)}^N,$$
where $Z$ is as above.
This is equivalent to
$$\langle Z,T \rangle = (1-q)^N.$$

Let $\gamma_1,\dots,\gamma_N \co S^1 \to D$ be
the figures of eight used to define $T$,
and let $E_1,\dots,E_N$ be the edges used to define $Z$.
Then $\gamma_i$ intersects $E_i$ at two points.
Call these points $y^+_i$ and $y^-_i$,
where $y^+_i$ is above $y^-_i$.
Then $T$ and $Z$ intersect at the $2^N$ points
$$(y^\pm_1,\dots,y^\pm_N).$$
Each such point $\mathbf{y}$ contributes a monomial
$\pm q^k$ to $\langle T,Z \rangle$.

We can assume that our basepoint of $T$ is given by
$$\mathbf{t} = (y^-_1,\dots,y^-_N).$$
The orientation of the intersection
of $E_i$ and $\gamma_i$ at $y^-_i$ is positive.
Thus $\mathbf{t}$ contributes $1$ to $\langle T,Z \rangle'$.

Now suppose $\mathbf{y}$ and $\mathbf{y}'$
are two points of intersection between $Z$ and $T$
that differ only at the mobile point of color $i$,
where $\mathbf{y}$ has $y^-_i$
and $\mathbf{y}'$ had $y^+_i$.
Let $\xi$ be the loop in $C$
that follows a path in $T$ from $\mathbf{y}$ to $\mathbf{y}'$,
and then follows a path in $Z$ back to $\mathbf{y}'$.
This can be represented by a braid
in which all strands are straight
except the strand of color $i$,
which makes a positive full twist
around the strands of color $i+1,\dots,N+1$.
Thus $\rho_\mcolor(\xi) = q$.
Also note that the orientation of the intersection at $\mathbf{y}'$
is the opposite of that at $\mathbf{y}$.
Thus
if $\mathbf{y}$ contributes $\pm q^k$ to $\langle T,Z \rangle'$
then $\mathbf{y}'$ contributes $\mp q^{k+1}$.

Summing the contributions of the $2^N$ points in $T \cap Z$ we obtain
$$\langle T,Z \rangle' = (1 - q)^N,$$
as required.
\end{proof}

\begin{lem}
\label{lem:TZ}
If $n=1$ and $Z$ is as above
then the image of $T$ in $H_m(C,\partial C)$ is $(1+q+\dots+q^N)Z$.
\end{lem}

\begin{proof}
We have $\langle S,Z \rangle = 1$.
Thus it suffices to prove the identity
$$\langle S,T \rangle = 1+q+\dots+q^N.$$

\begin{figure}
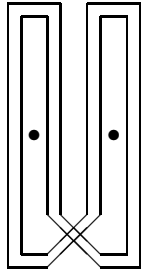

$$\barify$$
\caption{A stretched version of $T$.}
\label{fig:barify}
\end{figure}

Let $\gamma_1,\dots,\gamma_N$ be the figures of eight used to define $T$.
Isotope $T$ so that the disk $X$,
as defined in Section \ref{subsec:X},
is below the interval $[p_1,p_2]$.
See Figure \ref{fig:barify}.
Now each $\gamma_i$ intersects the interval $[p_1,p_2]$
at two points $a_i$ and $b_i$,
where $a_i$ is to the left of $b_i$.
Thus the points $a_1,\dots,a_N,b_1,\dots,b_N$
are in order from left to right.

For $i=0,1,\dots,N$,
let
$$\mathbf{y}_i = (a_1,\dots,a_i,b_{i+1},\dots,b_N).$$
Then $\mathbf{y}_0,\dots,\mathbf{y}_N$
are the points of intersection between $S$ and $T$.
Each of these contributes a monomial $\pm q^k$ to $\langle S,T \rangle$.

The sign of the intersection of $S$ and $T$ at $\mathbf{y}_i$ is positive
for all $i=0,\dots,N$.
It is not hard to see that
$\mathbf{y}_N$ contributes $+1$ to $\langle S,T \rangle$.
Let $\xi$ be a loop in $C$
that follows a path in $T$ from $\mathbf{y}_i$ to $\mathbf{y}_{i-1}$,
and then follows a path in $S$ back to $\mathbf{y}_i$.
This is the loop where
all mobile points remain stationary except for the point of color $i$,
which moves along $\gamma_i$ from $a_i$ to $b_i$,
and then horizontally back to $a_i$.
Then $\rho_\mcolor(\xi) = q$.
Thus if $\mathbf{y}_i$ contributes $q^k$ to $\langle S,T \rangle$
then $\mathbf{y}_{i-1}$ contributes $q^{k+1}$.
Summing the contributions of $\mathbf{y}_i$ for all $i$
gives the desired identity.
\end{proof}

\section{A partial barcode}
\label{sec:partialbarcode}

The aim of this section is
to prove a certain identity in $H_m^\lf(C,\partial C)$,
which will show that $S$ can be,
in some sense,
partially converted into a barcode.

Recall that $S$ was defined to be the product of
$n$ copies of an $N$-dimensional ball.
Call these $N$-balls $S_1,\dots,S_N$.
Let $Z$ be the nontrivial barcode for the case $n=1$,
as defined in Section \ref{subsec:images}.
Let $Z_i$ be the product
$$Z_i = S_1 \times \dots \times S_{i-1} \times
        Z \times S_{i+1} \times \dots \times S_N.$$
The basepoint $\mathbf{s}$ lies in $Z_i$,
so the path $\zeta$ determines a lift of $Z_i$ to $\Local$.
We obtain an element of $H_m^\lf(C,\partial C)$,
which we also call $Z_i$.
The aim of this section is to prove the following.

\begin{lem}
\label{lem:SZ}
$(q-1)^N S = (1+q+\dots+q^N) Z_i$ in $H_m^\lf(C,\partial C)$.
\end{lem}

First consider the case $n=1$.
By Lemma \ref{lem:TS},
$T = (q-1)^N S$.
It remains to show that
$$T = (1+q+\dots+q^N) Z_1.$$
But this is immediate from Lemma \ref{lem:TZ}.
We now describe how we could obtain this identity
in a way that will generalize to $n>1$.

First,
vertically ``stretch'' $T$,
as suggested by Figure \ref{fig:barify}.
Continue this stretching process and use an excision argument
to obtain a disjoint union of barcodes.
One of these must be $Z_1$,
with the desired coefficient.
Any other barcode must correspond to a trivial code sequence.
Such a barcode represents zero in $H_m^\lf(C,\partial C)$,
since one of the vertical edges can be slid to the boundary of the disk.

We can apply most of this argument to the case $n>1$.
The only difficulty is that the $N$-balls
$S_1,\dots,S_{i-1}$ and $S_{i+1},\dots,S_n$
prevent us from simply sliding a vertical edge to the boundary of the disk.
To overcome this problem,
we prove a claim that will imply that each such $N$-ball
is in some sense ``transparent'' to any other mobile point.
We need to make some definitions
before we can state the claim precisely.

Fix any $j=1,\dots,N$.
Let
$$\mcolor' = (1,2,\dots,N,j).$$
Let $C'$ be the configuration space
$$C' = C_{\mcolor'}(D_{(0,N+1)}).$$
Let $S'$ be the product of
the usual $N$-ball in $C_{(1,\dots,N)}(D_{(0,N+1)})$
and a circle of color $j$ around the interval $[p_1,p_2]$.
This is an $(N+1)$-dimensional submanifold of $C'$.

Let $g$ be the generator of $\pi_1(S')$.
Then $g$ can be represented by a braid
with strands of colors $0,1,\dots,N+1$ that are straight,
and a strand of color $i$ that
makes a positive full twist around all of the other strands.
Then
$$\rho_\mcolor (g) = (-q^{-1})^2(q^{1/2})^4 = 1.$$
Thus we can lift of $S'$ to $\Local$.
This represents an element of $H_m^\lf(C)$,
which we also call $S'$.

\begin{clm}
$S' = 0$ in $H_m^\lf(C')$.
\end{clm}

\begin{proof}
First consider the case $N=1$.
Then $S'$ is simply
the product of an edge $\gamma$ between the two puncture points
and a circle $\delta$ around $\gamma$.

Let $Z$ be the barcode corresponding to the code sequence $(0,1,1,2)$.
It suffices to show that
$$\langle S',Z \rangle' = 0.$$
Let $E$ and $E'$ be properly embedded vertical edges
passing between the puncture points,
where $E'$ is to the right of $E$.
As a closed $2$-ball,
$Z$ is the product of $E$ and $E'$,
both having color $1$.

\begin{figure}
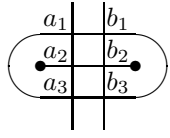

$$\transparent$$
\caption{$S'$ and $Z$ in the case $N=1$.}
\label{fig:transparent}
\end{figure}

Let $a_1$, $a_2$ and $a_3$ be the points of intersection
between $E$ and $\gamma \cup \delta$,
reading from top to bottom.
Let $b_1$, $b_2$, and $b_3$ be the analogous points of intersection
between $E'$ and $\gamma \cup \delta$.
See Figure \ref{fig:transparent}.
There are four points of intersection between $S'$ and $Z$,
namely
\begin{itemize}
\item $\mathbf{y}_1 = (a_2,b_1)$,
\item $\mathbf{y}_2 = (a_2,b_3)$.
\item $\mathbf{y}_3 = (a_1,b_2)$,
\item $\mathbf{y}_4 = (a_3,b_2)$,
\end{itemize}
Each of these contributes a monomial $\pm q^k$ to $\langle S',Z \rangle'$.
Assume the orientations and lifts to $\Local$
were chosen so that $\mathbf{y}_1$ contributes $+1$.

For $i \in \{1,2,3,4\}$,
let $\xi_i$ be a loop
that follows a path in $Z$ from $\mathbf{y}_1$ to $\mathbf{y}_i$,
and then follows a path in $S'$ back to $\mathbf{y}_1$.
Then
\begin{itemize}
\item $\rho_\mcolor(\xi_1) = 1$,
\item $\rho_\mcolor(\xi_2) = (q^{1/2})^2 = q$,
\item $\rho_\mcolor(\xi_3) = (-q^{-1})^{-1} = -q$,
\item $\rho_\mcolor(\xi_4) = (-q^{-1})(q^{1/2})^2 = -1$.
\end{itemize}

For $i = 1,2,3,4$,
let $\epsilon_i$ be the sign of the intersection
of $S'$ and $Z$ at $\mathbf{y}_i$.
By assumption,
$\epsilon_1 = 1$.
For $j = 1,2,3$,
the intersections of the relevant edges at $a_j$ and $b_j$
have the same sign.
The intersections at $a_1$ and $a_3$ have opposite signs.
The loops $\xi_3$ and $\xi_4$ transpose the two mobile points.
Combining these facts, we obtain
$\epsilon_4 = 1$ and $\epsilon_2 = \epsilon_3 = -1$.
Thus
$$\langle S',Z \rangle' = 1 - q - (-q) + (-1) = 0.$$

Now consider the case $N > 1$.
The only nontrivial code sequence is
$$(0,1,\dots,j-1,j,j,j+1,\dots,N,N+1).$$
Let $Z$ be the corresponding barcode.
We must show that
$$\langle S',Z \rangle' = 0.$$
As a closed $m$-ball,
$Z$ is a product of vertical edges
$$E_1,\dots,E_{j-1},E_j,E'_j,E_{j+1},\dots,E_N.$$
Here,
$E_k$ has color $k$ for $k=1,\dots,N$,
and $E'_j$ has color $j$.
Let $y_k$ be the point of intersection
between $E_k$ and the interval $[p_1,p_2]$.
Any point of intersection between $S'$ and $Z$
must include the mobile points $y_k$ of color $k$
for every $k \neq j$.
These points play no important role
since they remain the same throughout the proof.
The rest of the computation proceeds exactly as in the case $N=1$.

This completes the proof of the claim,
and hence of the lemma.
\end{proof}

\section{Bridge-preserving isotopy}
\label{sec:bridgepreserving}

We use the notation
$$\sigma_{2112} = \sigma_2 \sigma_1^2 \sigma_2.$$
The aim of this section is to prove the following.

\begin{lem}
\label{lem:bridgepreserving}
$Q(\sigma_{2112} \beta) = Q(\beta \sigma_{2112}) = Q(\beta)$.
\end{lem}

Combined with Lemma \ref{lem:heightpreserving},
this implies that $Q(\beta)$ is invariant under any isotopy of $\hat{\beta}$
through links that are in bridge position.

\begin{clm}
$Q(\sigma_{2112} \beta) = Q(\beta)$.
\end{clm}

\begin{proof}
We have
$$\rho_\pcolor(\sigma_{2112} \beta) = q^{-1} \rho_\pcolor(\beta).$$ 
By this and the properties of the sesquilinear pairing,
it suffices to show that the identity
$$\sigma_{2112} S = q S$$
holds in $H_m^\lf(C,\partial C)$.

Let $Z_2$ be as defined in Section \ref{sec:partialbarcode}.
By Lemma \ref{lem:SZ},
it suffices to show that $\sigma_{2112} Z_2 = q Z_2$.
We can choose the function $\sigma_{2112}$ to act as the identity on
the subset $Z_2$ of $C$.
It remains to show that $\sigma_{2112}$ acts as multiplication by $q$
on the fiber over $\mathbf{s}$.

Let $\xi$ be the concatenation of the paths
$\sigma_{2112} \zeta$ and $\overline{\zeta}$.
We can represent $\xi$ by a braid
in which strands of color $1,\dots,N$
wind in parallel around a strand of color $0$.
Thus $\rho_\mcolor(\xi) = q$.
Thus $\sigma_{2112} Z_2 = q Z_2$,
as required.
\end{proof}

It remains to show that $Q(\beta \sigma_{2112}) = Q(\beta)$.
It suffices to prove the following.

\begin{clm}
$Q(\beta^{-1}) = \overline{Q(\beta)}$.
\end{clm}

\begin{proof}
We have the following identities.
\begin{itemize}
\item $q^{m/2} = q^m \overline{(q^{m/2})}$,
\item $[N+1] = \overline{[N+1]}$,
\item $\rho_\pcolor(\beta^{-1}) = \overline{\rho_\pcolor(\beta)}$.
\end{itemize}
By the definition of $Q(\beta)$,
it remains to show that
$$\langle S,\beta^{-1}(T) \rangle=q^m\overline{\langle S,\beta(T) \rangle}.$$
By Lemma \ref{lem:TS} and the properties of the sesquilinear pairing,
this is equivalent to
$$\langle \beta(T),T \rangle = (-1)^m\overline{\langle T,\beta(T) \rangle}.$$
This follows from the symmetry property of the pairing.
\end{proof}

\section{Markov-Birman stabilization}
\label{sec:markov}

Let $\pcolor'$ be the $(2n+2)$-tuple $(0,N+1,0,N+1,\dots,0,N+1)$.
Let
$$\iota \co B_\pcolor \to B_{\pcolor'}$$
be the obvious inclusion map.
The {\em Markov-Birman stabilization} of $\beta$
is the braid
$$\beta' = (\sigma_{n+1}^{-1} \sigma_n \sigma_{n+1}) \iota(\beta).$$
The aim of this section is to prove the following.

\begin{lem}
\label{lem:markov}
If $\beta'$ is the Markov-Birman stabilization of $\beta$
then $Q(\beta') = Q(\beta)$.
\end{lem}

Combined with
Lemmas \ref{lem:heightpreserving} and \ref{lem:bridgepreserving},
this implies that $Q(\beta)$ is an invariant of
the oriented knot or link $\hat{\beta}$.

We make the following definitions.
\begin{itemize}
\item $m' = m + N$,
\item $\mcolor'$ is the $m'$-tuple $(1,\dots,N,1,\dots,N,\dots,1,\dots,N)$,
\item $D' = D_{\pcolor'}$,
\item $C' = C_{\mcolor'}(D')$,
\item $S'$ and $T'$ are the obvious
      embedded $m'$-ball and immersed $m'$-torus in $C'$.
\end{itemize}
We have the identities
$$\rho_\pcolor(\beta') = q^{N/2} (\rho_\pcolor(\beta)),$$
$$q^{m'/2} = q^{N/2} (q^{m/2}).$$
Thus it suffices to show
\begin{equation}
\label{eq:ma}
\langle S',\beta'(T') \rangle = \langle S,\beta(T) \rangle.
\end{equation}

Let $Z_n$ be the subset of $C$
as defined in Section \ref{sec:partialbarcode}.
Let $Z'_n$ be the subset of $C'$ defined similarly,
namely by replacing the second to rightmost $N$-ball of $S'$
by a barcode.
By Lemma \ref{lem:SZ},
equation (\ref{eq:ma}) is equivalent to
$$\langle Z'_n,\beta'(T') \rangle = \langle Z_n,\beta(T) \rangle.$$
This,
in turn,
is equivalent to
\begin{equation}
\label{eq:mb}
\langle \sigma(Z'_n),\iota(\beta)(T') \rangle = \langle Z_n,\beta(T) \rangle,
\end{equation}
where
$$\sigma = \sigma_{n+1}^{-1} \sigma_n^{-1} \sigma_{n+1}.$$

First let us look at $\sigma(Z'_n)$.
Let $D_3$ be the three times punctured disk
consisting of points in $D'$ to the right of
a vertical line between $p_{2n-1}$ and $p_{2n}$.
There is an embedding
$$C_\mcolor(D' \setminus D_3) \times C_{(1,\dots,N)}(D_3) \to C'.$$
We can assume that $Z_n$ lies in $C_\mcolor(D' \setminus D_3)$,
and
$$Z'_n = Z_n \times S_{n+1},$$
where $S_{n+1}$ is the obvious $N$-ball in $C_{(1,\dots,N)}(D_3)$.
Thus
$$\sigma(Z'_n) = Z_n \times \sigma(S_{n+1}).$$

Next we look at $\iota(\beta)(T')$.
Let $D_2$ be the twice punctured disk
consisting of points in $D'$ to the right of
a vertical line between $p_{2n}$ and $p_{2n+1}$.
There is an embedding
$$C_\mcolor(D' \setminus D_2) \times C_{(1,\dots,N)}(D_2) \to C'.$$
We can assume that $T$ lies in $C_\mcolor(D' \setminus D_2)$,
and
$$T' = T \times T_{n+1},$$
where $T_{n+1}$ is
the obvious $N$-dimensional torus in $C_{(1,\dots,N)}(D_2)$.
Then
$$\iota(\beta)(T') = \beta(T) \times T_{n+1}.$$

Equation (\ref{eq:mb}) is now equivalent to
$$ \langle Z \times \sigma(S_{n+1}), \beta(T) \times T_{n+1} \rangle
  = \langle Z,\beta(T) \rangle. $$
Any point of intersection between
$Z \times \sigma(S_{n+1})$ and $\beta(T) \times T_{n+1}$
must lie in
$$C_\mcolor(D' \setminus D_3) \times C_{(1,\dots,N)}(D_2),$$
which is the intersection of the two relevant product spaces.
Thus it suffices to show
\begin{equation}
\label{eq:mc}
\langle \sigma(S_{n+1}), T_{n+1} \rangle = 1.
\end{equation}
We can take this intersection pairing to be between
submanifolds of $C_{(1,\dots,N)}(D_3)$ .

\begin{figure}
$$\markov$$
\caption{Computing $\langle \sigma(S_{n+1}), T_{n+1} \rangle$ when $N=1$.}
\label{fig:markov}
\end{figure}

Equation (\ref{eq:mc}) follows from
a direct computation of a particular intersection pairing.
Figure \ref{fig:markov} shows the case $N=1$.
The case $N > 1$ is similar.
There is one point of intersection
$\mathbf{y}$ between $\sigma(S_{n+1})$ and $T_{n+1}$.
The sign of this intersection is positive.
Both $\sigma(S_{n+1})$ and $T_{n+1}$ come with an path
from a configuration of points on $\partial D_3$ to $\mathbf{y}$.
These paths are homotopic relative to endpoints.
This completes the proof of equation (\ref{eq:mc}),
and hence of the lemma.

\section{The skein relation}
\label{sec:skein}

Let $\beta_+ = \sigma_2^{-1} \sigma_1 \sigma_2 \beta$
and $\beta_- = \sigma_2^{-1} \sigma_1^{-1} \sigma_2 \beta$.
The aim of this section is to prove the following.

\begin{lem}
\label{lem:skein}
$q^{(N+1)/2} Q(\beta_-) - q^{-(N+1)/2} Q(\beta_+)
 = (q^{1/2} - q^{-1/2}) Q(\beta)$.
\end{lem}

We have the identities
\begin{eqnarray*}
\rho_\pcolor(\beta_+) &=& q^{N/2} \rho_\pcolor(\beta), \\
\rho_\pcolor(\beta_-) &=& q^{-N/2} \rho_\pcolor(\beta).
\end{eqnarray*}
Thus it suffices to show
$$ q^{1/2} \langle S,\beta_-(T) \rangle
   - q^{-1/2}\langle S,\beta_+(T) \rangle
   = (q^{1/2} - q^{-1/2}) \langle S,\beta(T) \rangle.  $$
Let $Z_2$ be as defined in Section \ref{sec:partialbarcode}.
By Lemma \ref{lem:SZ},
it suffices to show that
$$ q^{1/2} \langle Z_2,\beta_-(T) \rangle
   - q^{-1/2}\langle Z_2,\beta_+(T) \rangle
   = (q^{1/2} - q^{-1/2}) \langle Z_2,\beta(T) \rangle.  $$
By some simple manipulation, this is equivalent to
$$ \langle \sigma_2^{-1} (\sigma_1 - 1) (1 + q \sigma_1^{-1}) \sigma_2 (Z_2),
        \beta(T) \rangle = 0.  $$
Thus it suffices to prove the identity
\begin{equation}
\label{eq:sa}
\sigma_2^{-1} (\sigma_1 - 1) (1 + q \sigma_1^{-1}) \sigma_2 (Z_2) = 0
\end{equation}
in $H_m^\lf(C,\partial C)$.

Let $D_3$ be the set of points in $D$
on or to the left of a vertical line between $p_3$ and $p_4$.
Let
$$C_1 = C_{(1,\dots,N)}(D_3).$$
Let $\mcolor_2$ be the $(m-N)$-tuple
$$\mcolor_2 = (1,\dots,N,1,\dots,N,\dots,1,\dots,N),$$
and let
$$C_2 = C_{\mcolor_2}(D \setminus D_3).$$
There is an obvious embedding
$$C_1 \times C_2 \to C.$$
We can write
$$Z_2 = S_1 \times Z',$$
where $S_1$ is the obvious $N$-ball in $C_1$,
and $Z'$ is an $(N-m)$-manifold in $C_2$.

Now $\sigma_1$ and $\sigma_2$ both act as the identity on $D \setminus D_3$.
Thus,
to prove equation (\ref{eq:sa}),
it suffices to show that
\begin{equation}
\label{eq:sb}
\sigma_2^{-1} (\sigma_1 - 1)(1 + q \sigma_1^{-1}) \sigma_2 (S_1) = 0
\end{equation}
in $H_N^\lf(C_1)$.

We now eliminate the conjugation by $\sigma_2$ in equation (\ref{eq:sb}).
Let $D'_3 = \sigma_2 D_3$.
This is a disk with three puncture points,
which have colors $0,0,N+1$,
reading from left to right.
Let
$$C'_1 = C_{(1,\dots,N)} D'_3.$$
Let $S'_1 = \sigma_2 S_1$.
Then equation (\ref{eq:sb}) is equivalent to the identity
\begin{equation}
\label{eq:sc}
(\sigma_1 - 1)(1 + q \sigma_1^{-1}) (S'_1) = 0
\end{equation}
in $H_N^\lf(C'_1)$.

In this setting,
there are only two nontrivial code sequence,
namely
$(0,1,2,\dots,N+1,0)$ and $(0,N+1,N,\dots,1,0)$.

Suppose $Z$ is the barcode corresponding to $(0,1,2,\dots,N+1,0)$.
Then $S'_1$ and $Z$ do not intersect,
so
$$\langle S'_1,Z \rangle' = 0.$$
Now $\sigma_1(S'_1)$ and $Z$ intersect at a single point $\mathbf{y}$.
Similarly,
$\sigma_1^{-1}(S'_1)$ and $Z$ intersect at a single point,
which we can assume is also $\mathbf{y}$.
The signs of these intersections are the same.
Each of $\sigma_1(S'_1)$ and $\sigma_1^{-1}(S'_1)$ comes with a path
from $\mathbf{y}$ to $\mathbf{x}$.
These paths differ by
the direction the points of colors $1,\dots,N$
pass around the middle puncture point.
Thus
$$\langle\sigma_1^{-1}(S'_1),Z\rangle' = q\langle\sigma_1(S'_1),Z\rangle'.$$
A simple computation now gives
$$\langle (\sigma_1 - 1)(1 + q \sigma_1^{-1}) (S'_1), Z \rangle' = 0.$$

Now suppose $Z$ is the barcode corresponding to $(0,N+1,N,\dots,1,0)$.
Then $\sigma_1$ acts as the identity on $Z$.
It follows that
$$\langle S'_1,Z \rangle'
  = \langle \sigma_1(S'_1),Z \rangle'
  = \langle \sigma_1^{-1}(S'_1),Z \rangle'.$$
A simple computation now gives
$$\langle (\sigma_1 - 1)(1 + q \sigma_1^{-1}) (S'_1),Z \rangle' = 0.$$

This completes the proof of equation (\ref{eq:sc}),
and hence of the lemma.

\section{Conclusion}
\label{sec:conclusion}

We are now ready to prove the main theorem of this paper.
Let $Q(\beta)$ be as defined in Section \ref{sec:def}.
Let $\hat{\beta}$ be the plat closure of $\beta$,
as an oriented knot or link.
Let $P(\hat{\beta})$ be the invariant of $\hat{\beta}$ of type $A_N$,
as defined in the introduction.

\begin{thm}
$Q(\beta) = P(\hat{\beta})$.
\end{thm}

\begin{proof}
Birman \cite{jB76} proved that
two braids have isotopic plat closures
if and only if they are related by a sequence of moves of certain types.
The original theorem applied to unoriented knots,
whereas we wish to apply it to oriented knots and links.
However the result is essentially the same.
The moves are those given in Lemmas \ref{lem:heightpreserving},
\ref{lem:bridgepreserving}, and \ref{lem:markov}.
Thus $Q(\beta)$ is an invariant of the oriented knot or link $\hat{\beta}$.

Suppose we have three links as shown in the skein relation
given at the beginning of this paper.
By applying an isotopy,
we can present these links as the plat closures of braids
$\beta_+$, $\beta_-$ and $\beta$ as in Section \ref{sec:skein}.
By Lemma \ref{lem:skein},
the invariant $Q$ satisfies the required skein relation.

It remains only to prove
that $Q$ is correctly normalized to take the value one for the unknot.
Suppose $n=1$ and $\beta$ is the identity braid.
By Lemma \ref{lem:TZ},
$$\langle S,T \rangle = 1+q+ \dots + q^N.$$
Thus
$$Q(\beta) = \frac{1}{[N+1]q^{N/2}}(1+q+ \dots + q^N) = 1,$$
as required.
\end{proof}

We now show how to
eliminate the factor of $[N+1]$ from the definition of $Q$.
Suppose the rightmost strand of $\beta$
makes no crossings with any other strands.
Note that
any oriented knot or link is the plat closure of some such braid $\beta$.
Let $m' = m - N$.
Let $\mcolor'$ be the $m'$-tuple
$$\mcolor' = (1,\dots,N,1,\dots,N,\dots,1,\dots,N).$$
Let $C' = C_{\mcolor'}(D_\pcolor)$.
Recall that $S$ was defined to be the product of
$n$ copies of an $N$-dimensional ball.
Let $S'$ be the product of all but the rightmost of these $N$-balls,
as a subset of $C'$.
Similarly,
let $T'$ be the product of
all but the rightmost $N$-torus used to define $T$.
Let
$$Q'(\beta) = \rho_\pcolor(\beta) q^{-m'/2} \langle S',\beta(T') \rangle.$$

\begin{thm}
If the rightmost strand of $\beta$ makes no crossings with any other strands
then $Q'(\beta) = P(\hat{\beta})$.
\end{thm}

\begin{proof}
Let $Z_n$ be as defined in Section \ref{sec:partialbarcode}.
By Lemma \ref{lem:SZ},
$$\overline{(q-1)}^N \langle S,\beta(T) \rangle
  = \overline{(1+q+\dots+q^N)} \langle Z_n,\beta(T) \rangle.  $$

Let $D_1$ be the once punctured disk
consisting of points in $D$ to the right of
a vertical line between $p_{2n-1}$ and $p_{2n}$.
There is an embedding
$$C_{\mcolor'}(D \setminus D_1) \times C_{(1,\dots,N)}(D_1) \to C.$$
Then $Z_n = S' \times Z$,
where $Z$ is the obvious barcode in $C_{(1,\dots,N)}(D_1)$.

Let $D_2$ be the twice punctured disk
consisting of points in $D$ to the right of
a vertical line between $p_{2n-2}$ and $p_{2n-1}$.
There is an embedding
$$C_{\mcolor'}(D \setminus D_2) \times C_{(1,\dots,N)}(D_2) \to C.$$
Then $T = T' \times T_n$,
where $T_n$ is the obvious $N$-torus in $C_{(1,\dots,N)}(D_2)$.

By assumption,
$\beta$ acts as the identity on $D_1$.
Any point of intersection between $Z_n$ and $\beta(T)$
must lie in
$$C_{\mcolor'}(D \setminus D_2) \times C_{(1,\dots,N)}(D_1),$$
which is the intersection of the two relevant product spaces.
Thus
$$\langle Z_n,\beta(T) \rangle
   = \langle S',\beta(T') \rangle  \langle Z,T_n \rangle. $$
By Lemma \ref{lem:TS},
$$\langle Z,T_n \rangle = (1-q)^N.$$
A straightforward calculation now gives $Q'(\beta) = Q(\beta)$,
as required.
\end{proof}

The computational definition of the pairing works
over any ring containing an invertible element $q$.
Thus $Q'(\beta)$ is well defined
over any ring containing an invertible element $q^{1/2}$.
Since it is a polynomial in $q^{\pm 1/2}$,
the above theorem applies for any such ring.


\newcommand{\etalchar}[1]{$^{#1}$}
\providecommand{\bysame}{\leavevmode\hbox to3em{\hrulefill}\thinspace}
\providecommand{\MR}{\relax\ifhmode\unskip\space\fi MR }
\providecommand{\MRhref}[2]{%
  \href{http://www.ams.org/mathscinet-getitem?mr=#1}{#2}
}
\providecommand{\href}[2]{#2}

\end{document}